\documentclass[a4paper,12pt,reqno]{amsart}
\usepackage{amssymb}
\usepackage[dvips]{graphicx}
\usepackage{hyperref}

\input epsf.tex
\newdimen\xsize
\newdimen\oldbaselineskip
\newdimen\oldlineskiplimit
\def\restorelineskip{\baselineskip=\oldbaselineskip%
\lineskiplimit=\oldlineskiplimit}
\def\putm[#1][#2]#3{
\hbox{\vbox to 0pt{\parindent=0pt%
\vskip#2\xsize\hbox to0pt{\hskip#1\xsize $#3$\hss}\vss}}}%
\long\def\Line#1{\hbox to \hsize{#1}}
\def\putt[#1][#2]#3{
\vbox to 0pt{\noindent\hskip#1\xsize\lower#2\xsize%
\vtop{\restorelineskip#3}\vss}}

\makeatletter
\def\xbig[#1]#2{{\hbox{$\m@th\left#2\vbox to#1\xsize{}%
\right.\n@space$}}}
\makeatother
\def\xlar[#1]#2{%
\smash{\mathop{ \hbox to #1\xsize{\leftarrowfill}}\limits^{#2}}}
\def\xrar[#1]#2{%
\smash{\mathop{ \hbox to #1\xsize{\rightarrowfill}}\limits^{#2}}}
\def\xline[#1]{\hbox to #1\xsize{\leaders\hrule\hfill}}

\DeclareFontFamily{U}{rsf}{\skewchar\font'177}%
\DeclareFontShape{U}{rsf}{m}{n}{<-6>rsfs5<6-8>rsfs7<8->rsfs10}{}%
\DeclareFontShape{U}{rsf}{b}{n}{<-6>rsfs5<6-8>rsfs7<8->rsfs10}{}%
\DeclareMathAlphabet\RSFS{U}{rsf}{m}{n}
\SetMathAlphabet\RSFS{bold}{U}{rsf}{b}{n}
\def\sf#1{{\mathsf{#1}}}

\def\slsf{\slshape \sffamily }

\def\msmall#1{\mathchoice{\hbox{\small$\displaystyle {#1}$}}{#1}{#1}{#1}}

\def\cc{{\mathbb C}}

\def\rr{{\mathbb R}}

\def\sph{{\mathbb S}}
\def\nn{{\mathbb N}}

\def\pp{{\mathbb P}}

\def\arg{{\sf{Arg}}}

\def\area{\sf{area}}

\def\hol{\sf{Hol}}

\def\id{\sf{Id}}

\def\inter{\sf{Int}}

\def\ker{\sf{Ker}\,}
\def\lim{\mathop{\sf{lim}}}

\def\max{\sf{max}}

\def\ord{\sf{ord}}

\def\var{{\sf{Var}}}

\def\eps{\varepsilon}

\def\<{\langle}\let\la=\<
\def\>{\rangle}\let\ra=\>
  
\def\comp{\Subset}
\def\d{\partial}

\def\ddef{\mathrel{{=}\raise0.3pt\hbox{:}}}
\def\deff{\mathrel{\raise0.3pt\hbox{\rm:}{=}}}

\def\fraction#1/#2{\mathchoice{{\msmall{ #1\over#2}}}%
{{ #1\over #2 }}{{#1/#2}}{{#1/#2}}}
\def\norm#1{\left\Vert{#1}\right\Vert}

\def\le{\leqslant}

\def\emptyset{\varnothing}

\def\longpoints{\leaders\hbox to 0.5em{\hss.\hss}\hfill \hskip0pt}
\def\stateskip{\smallskip}
\def\state#1. {\stateskip\noindent{\bf#1. }} 
\def\statep#1. {\stateskip\noindent{\bf#1 }} 
\def\proof{\state Proof. \2}

\def\Chi{\raise 2pt\hbox{$\chi$}}

\def\ie{\hskip1pt plus1pt{\sl i.e.\/,\ \hskip1pt plus1pt}}

\def\sli{{\sl i)} } 
\def\slii{{\sl i$\!$i)} } 
\def\sliii{{\sl i$\!$i$\!$i)} }

\def\Chi{\raise 2pt\hbox{$\chi$}}

\let\phI=\phi\let\phi=\varphi\let\varphi=\phI
\let\cal=\mathcal
\def\cala{{\cal A}}
\def\calb{{\cal B}}
\def\calc{{\cal C}}

\def\calf{{\cal F}}

\def\calp{{\cal P}}

\def\calr{{\cal R}}

\def\calx{{\cal X}}

\def\calz{{\cal Z}}

\def\eps{\varepsilon}

\def\comp{\Subset}
\def\d{\partial}

\def\1{{1\mkern-5mu{\rom l}}}

\def\ge{\geqslant}

\def\fraction#1/#2{\mathchoice{{\msmall{ #1\over#2}}}%
{{ #1\over #2 }}{{#1/#2}}{{#1/#2}}}

\def\le{\leqslant}

\def\emptyset{\varnothing}
\newcommand{\2}{\thinspace}

\def\qed{\ \ \hfill\hbox to .1pt{}\hfill\hbox to .1pt{}\hfill $\square$\par}

\def\comment#1\endcomment{}

 \belowdisplayskip=\abovedisplayskip

\def\lineeqqno(#1){\hfill\llap{\vbox to 10pt%
{\vss\begin{align} \eqqno(#1)\end{align}\vss}}\vskip1pt}

\advance\headheight 1.2pt

\def\ShowwLLabel#1{}

\def\thechpt{\Roman{chpt}}

\def\newchapt[#1]#2{\newpage%
\refstepcounter{chpt}\setcounter{subsection}{0}%
\setcounter{thm}{0}\setcounter{defi}{0}%
\setcounter{rema}{0}\setcounter{exrc}{0}%
\renewcommand{\thesubsection}{\thechpt.\arabic{subsection}}%
\section*{\begin{center}\huge \bf Chapter \thechpt\\
#2 \end{center}}\label{#1}%
\ \smallskip%
\markboth{Chapter \thechpt}{#2}%
}
\def\newsect[#1]#2{\refstepcounter{section}\setcounter{equation}{0}%
\renewcommand{\thesubsection}{\arabic{section}.\arabic{subsection}}%
\section*{\arabic{section}.
#2}\vspace{-20pt}\label{#1}\vspace{20pt}%
\markboth{Section \arabic{section}}{#2}}

\def\newlect[#1]#2{\refstepcounter{section}%
\renewcommand{\thesubsection}{\arabic{section}.\arabic{subsection}}%
\section*{Lecture \arabic{section}\\
#2}\label{#1}%
\markboth{Lecture \arabic{section}}{#2}}

\def\newprg[#1]#2{\refstepcounter{subsection}%
\subsection*{{\thesubsection.\ #2}} \label{#1}%
}

\def\newappx[#1]#2{%
\refstepcounter{appx}\setcounter{section}{0}%
\renewcommand{\thesubsection}{A\arabic{appx}.\arabic{subsection}}%
\section*{Appendix \arabic{appx}\\ #2}
\label{#1}%
\markboth{Appendix A\arabic{appx}}{#2}
}

\newtheorem{thm}{Theorem}[section]
   \def\newthm#1{\begin{thm}\label{#1}}

\newtheorem{nnthm}{Theorem}
   \def\newthm#1{\begin{nnthm}\label{#1}}

\newtheorem{lem}{Lemma}[section]
   \def\newlemma#1{\begin{lem} \label{#1}}

\newtheorem{prop}{Proposition}[section]
   \def\newprop#1{\begin{prop}\label{#1}}

\newtheorem{nnprop}{Proposition}
   \def\newprop#1{\begin{nnprop}\label{#1}}

\newtheorem{corol}{Corollary}[section]
   \def\newcorol#1{\begin{corol} \label{#1}}

\newtheorem{nncorol}{Corollary}
   \def\newcorol#1{\begin{nncorol} \label{#1}}

\newtheorem{defi}{Definition}[section]
   \def\newdefi#1{\begin{defi} \label{#1}\rm }

\newtheorem{exmp}{Example}[section]
   \def\newexmp#1{\begin{exmp} \label{#1}\rm }

\newtheorem{nnexmp}{Example}
   \def\newexmp#1{\begin{nnexmp} \label{#1}\rm }

\newtheorem{exrc}{Exercise}
   \def\newexrc#1{\begin{exrc} \label{#1}\rm }

\newtheorem{rema}{Remark}[section]
   \def\newrema#1{\begin{rema} \label{#1}\rm }

\newtheorem{nnrema}{Remark}
   \def\newrema#1{\begin{nnrema} \label{#1}\rm }

\def\eqqno(#1){\label{(#1)}}
\def\eqqref(#1){(\ref{(#1)})}

\title{Banach Analytic Sets and a Non-Linear Version\\ of the Levi Extension Theorem}
\author{S. Ivashkovich}
\date{\today}

\address{
Universit\'e de Lille-1, UFR de Math\'ematiques, 59655 Villeneuve
d'Ascq, France} \email{ivachkov@math.univ-lille1.fr}
\address{IAPMM Nat. Acad. Sci. Ukraine
Lviv, Naukova 3b, 79601 Ukraine}

\keywords{Meromorphic function, analytic disk, Levi's theorem, Banach
analytic set.}
\thanks{This research is supported in part by the grant ANR-10-BLAN-0118}
\subjclass[2000]{Primary - 32D99, Secondary - 32H04}

\begin{document}

\begin{abstract}
We prove a certain non-linear version of the Levi extension theorem
for  meromorphic functions. This means that the meromorphic function
in question is supposed to be extendable along a sequence of complex
curves, which are arbitrary, not necessarily straight lines.
Moreover, these curves are not supposed to belong to any finite
dimensional analytic family. The conclusion of our theorem is that
nevertheless the function in question meromorphically extends along
an (infinite dimensional) analytic family of complex curves and its
domain of existence is a pinched domain filled in by this analytic
family.
\end{abstract}

\maketitle

\newsect[INT]{Introduction}

\newprg[INT.res]{Statement of the main result}
By $(\lambda , z)$ we denote the standard coordinates in $\cc^{2}$.
For $\eps >0$ consider the following ring domain
\begin{equation}
\eqqno(ring)
R_{1+\eps }=\{ (\lambda ,z)\in \cc^2: 1 -\eps <|\lambda |<1 + \eps, |z| < 1\} =
A_{1-\eps, 1+\eps} \times \Delta ,
\end{equation}
\ie $R_{1+\eps}$ is the product of the annulus $A_{1-\eps,1 +
\eps}:=\{z\in\cc : 1 -\eps <|\lambda|<1+\eps\}$ with the unit disk
$\Delta$. Let a sequence of holomorphic functions $\{ \phi_{k}:
\Delta_{1+\eps} \to \Delta\}_{k=1}^{\infty}$ be given such that
$\phi_k$ converge uniformly on $\Delta_{1+\eps}$ to some $\phi_0 :
\Delta_{1 + \eps} \to \Delta$. We say that such sequence is a {\slsf
test sequence} if $(\phi_k -\phi_0)|_{\d \Delta}$ doesn't vanish for
$k>>0$ and
\begin{equation}
\eqqno(test) \var\arg_{\d\Delta}(\phi_k -\phi_0) \qquad\text{stays
bounded when} \qquad k\to + \infty .
\end{equation}
Denote by $C_{k}$ the graph of $\phi_{k}$ in
$\Delta_{1+\eps}\times\Delta$, by $C_0$ the graph of $\phi_0$.

\begin{nnthm}
\label{levi-nonlin} Let $f$ be a meromorphic function on $R_{1+\eps}$ and
$\{\phi_k\}_{k=1}^{\infty}$ a test sequence such that for every $k$ the restriction
$f|_{C_{k}\cap R_{1+\eps}}$ is well defined and extends to a meromorphic function on the
curve $C_{k}$ and that the number of poles counting with
multiplicities of these extensions is uniformly bounded. Then there
exists an analytic family of holomorphic graphs
$\{C_{\alpha}\}_{\alpha\in \cala}$ parameterized by a Banach ball
$\cala$ of infinite dimension such that:

\smallskip\sli $f|_{C_{\alpha}\cap R_{1+\eps}}$ extends to a meromorphic
function on $C_{\alpha}$ for every $\alpha\in \cala$ and the number

\quad of poles counting with multiplicities of these extensions is
uniformly bounded.

\smallskip\slii  Moreover $f$ meromorphically extends as a function of
two variables $(\lambda , z)$ to the

\smallskip\quad pinched domain  $\calp \deff \inter\left(\bigcup_{\alpha\in
\cala}C_{\alpha}\right)$ swept by $C_{\alpha}$.
\end{nnthm}

\begin{figure}[h]
\centering
\includegraphics[width=1.5in]{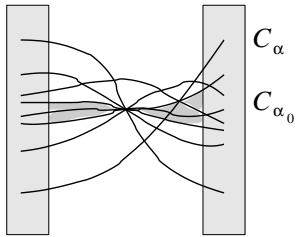}
\caption{The brighter dashed zone on this picture represents the
ring domain $R_{1+\eps}$ and curves are the graphs $C_{\alpha}$.
Around $C_{\alpha_0}$, the graph of $\phi_0 = \phi_{\alpha_0}$,
the analytic family $\{C_{\alpha}\}_{\alpha\in \cala}$ fills in an
another (darker) dashed zone, a pinched domain $\calp$. On this
picture there is exactly one pinch, the point at which most of
graphs intersect.} \label{pic1}
\end{figure}

Here by $C_{\alpha}$ we denote the graph of the function $\phi_{\alpha}$.
The notion of a {\slsf pinched domain}, though intuitively clear,
see Figure \ref{pic1}, is discussed in details at the beginning of
section \ref{PINCH}.

\begin{defi}
\label{gen-pos1}
Let's say that the graphs $\{ C_{k}\}$ of our functions $\{\phi_k\}$  are in {\slsf general position}
if for every point $\lambda_0\in \Delta$ there exists a  subsequence $\{\phi_{k_p}\}$ such that zeroes of
$\phi_{k_p}-\phi_0$ do not accumulate to $\lambda_0$.
\end{defi}

Theorem \ref{levi-nonlin}  implies the following {\slsf non-linear} Levi-type extension theorem:

\begin{nncorol}
\label{levi-c} If under the conditions of Theorem \ref{levi-nonlin}  curves
$\{C_{k}\}$ are in general position then $f$ extends to a
meromorphic function in the bidisk $\Delta_{1+\eps}\times \Delta$.
\end{nncorol}

\begin{nnrema} \rm
\label{rema-1} Let us explain the condition of a
general position. Take the sequence $C_{k}=\{ z=\frac{1}{k}\lambda
\}$ in $\cc^2$. Then the function $f(\lambda
,z)=e^{\frac{z}{\lambda}}$ is holomorphic in $R \deff \cc^*\times \cc$ and
extends holomorphically along every curve $C_{k}$. But it is not
holomorphic (even not meromorphic) in $\cc^2$. It is also
holomorphic when restricted to any curve $C= \{ z= \phi (\lambda)\}$
provided $\phi (0) =0$. Therefore the subspace $H_0$ of $\phi \in \hol
(\Delta_{1+\eps},\Delta)$ such that $f$ extends along the
corresponding curve is of codimension one. In fact this is the
general case: the Banach ball $\cala$ in Theorem \ref{levi-nonlin} appears
as a neighborhood of the limit point $\alpha_0$ in the subspace of
finite codimension of a well chosen Banach space of holomorphic
functions.
\end{nnrema}

\begin{nnrema} \rm 
\label{rema-2}
To explain the condition \eqqref(test) we impose on our test
sequences we construct in section \ref{EXMP} for the following
{\slsf non test} sequence $\phi_k(\lambda) = \left(\frac{2}{3}
\lambda \right)^k$ a holomorphic function $f$ in $\cc^*\times \cc$ 
which holomorphically extends along every graph $C_k$ but which is not
extendable meromorphically along any one-parameter analytic family
$\{\phi_{\alpha}\}$, see Example \ref{example1} there. 
\end{nnrema}

\newprg[INT.map]{Meromorphic mappings}
The assumption that $f$ is a function in Theorem \ref{levi-nonlin}
is not really important. We prove also a non-linear version of
an extension theorem for meromorphic mappings with values in general
complex spaces putting it into a form suitable for applications.
Let's call a family $\{\phi_t\in \hol (\Delta_{1+\eps},\Delta): t\in
T\}$ a {\slsf test family} if there exists $N\in \nn$ such that for
every pair $s\not=t\in T$ there exists a radius $1-\eps/2
<r<1+\eps/2$ such that $(\phi_s-\phi_t)|_{\d\Delta_r}$ doesn't
vanish and has winding number $\le N$. As usual by $C_t$ we denote
the graph of $\phi_{t}$.

\begin{nncorol}
\label{levi-gen} Let $X$ be a reduced, disk-convex complex
space and $f:R_{1+\eps}\to X$ a meromorphic mapping. Suppose that
there exists an uncountable test family of holomorphic functions
$\{\phi_{t}:\hol (\Delta_{1+\eps},\Delta ): t\in T\}$ such that
$f|_{C_{t}\cap R_{1+\eps}}$ holomorphically extends to $C_{t}$ for
every $t\in T$. Then $f$ extends to a meromorphic mapping from a
pinched domain $\calp$ to $X$.
\end{nncorol}
Moreover, there exists, like in Theorem \ref{levi-nonlin}, an
infinite dimensional family of graphs $C_{\alpha}$ parameterized by
a Banach ball $\cala$ such that $f|_{C_{\alpha}\cap R_{1+\eps}}$
holomorphically extends to $C_{\alpha}$ for all $\alpha\in\cala$.
The condition that $f|_{C_{t}\cap R_{1+\eps}}$ is assumed to extend
{\slsf holomorphically} should not be confusing because meromorphic
functions on curves are precisely holomorphic mappings to the
Riemann sphere $\pp^1$. I.e., the meromorphic functions case is the
case $X=\pp^1$ in this Corollary.

\newprg[INT.notes]{Structure of the paper, notes, acknowledgement.} 
{\slsf 1.} Theorem \ref{levi-nonlin} is proved in section \ref{PINCH}.  
The set $\cala$ such that $f|_{C_{\alpha}\cap R_{1+\eps}}$ meromorphically 
extends to $C_{\alpha}$ for $\alpha \in \cala$ is always a Banach analytic
subset of a neighborhood of $\phi_0$ in the Banach space $\hol
(\Delta_{1+\eps} ,\Delta)$. In particular the sequence $\{\phi_k\}$ of 
Example \ref{example1} is a Banach analytic set, namely the zero set of an 
appropriate singular integral operator, see section \ref{N-TEST} for more 
details. The (known) problem however is that a {\slsf every}
metrisable compact can be endoved with a structure of a Banach analytic set
(in an appropriate complex Banach space), see \cite{Mz}. For the case of a
converging sequence of points see Remark \ref{tot-real} for a very simple 
example. Therefore in infinite dimensional case from the
fact that our Banach analytic set contains a non-isolated point we cannot
deduce that it contains an analytic disk. Our major task here is to 
overcome this difficulty.

\smallskip\noindent{\slsf 2.} \rm In the case when $C_{k} = \Delta_{1+\eps}
\times \{z_k\}$ with $z_k\to 0$, \ie when $C_k$ are
horizontal disks this result is exactly the theorem of E. Levi, see
\cite{Lv} (case of holomorphic extension is due to Hartogs, see
\cite{Ht}). It should be said that Levi's theorem is usually stated
in the form as our Corollary \ref{levi-gen}: {\sl if $f$ as above extends
along an {\slsf uncountable family} of horizontal disks $C_{t} =
\Delta_{1+\eps}\times \{t\}$, then $f$ meromorphically extends
to $\Delta^2$.} But the proof goes as follows: one remarks that then
there exists a sequence $\{t_k\}$ (in fact an uncountable
subfamily) such that extensions along $C_{t_k}$-s have
uniformly bounded number of poles, and then the statement like our
Theorem \ref{levi-nonlin} is proved.

\smallskip\noindent{\slsf 3.} \rm If the number of poles of extensions
$f|_{C_k}$ is not uniformly bounded then the conclusion of
Theorem \ref{levi-nonlin} fails to be true even in the
case of horizontal disks. This is shown by the Example \ref{example2} 
in section \ref{EXMP}.

\smallskip\noindent{\slsf  4.} \rm In the case when $\{C_{t}\}_{t\in T}$ are
non-horizontal straight disks, \ie intersections of lines with
$\Delta^2$, Corollary \ref{levi-c} is due to Dinh, see \cite{Dh}
Corollaire 1. The proof in \cite{Dh} uses  results on the complex
Plateau problem in projective space (after an appropriate Segre
imbedding) and is essentially equivalent to the solution of this
problem. From the point of view of this paper this is a special
case when $\{C_k\}$ {\slsf ad hoc} belong to a finite dimensional
analytic family: in Levi case the family is one-dimensional, in the case
of Dinh two-dimensional. In section \ref{N-TEST}, after recalling the
necessary facts about singular integral transforms, we give a very
short proof of a non-linear extension theorem, see Theorem \ref{fin-dim-case},
in the case when $\{C_k\}$ are ad hoc included in an {\slsf arbitrary}
finite dimensional family.
In the straight case, \ie when  $\{C_{t}\}$ are
non-horizontal straight disks, the result of Corollary \ref{levi-gen}
for K\"ahler $X$ was proved in \cite{Sk} following the approach of \cite{Dh}.

\smallskip\noindent{\slsf 5.} \rm  It is important to outline that we
do not suppose {\slsf a priori} that $\{C_{k}\}$ are included into
any finite dimensional family of complex curves (ex. any family of
algebraic curves of uniformly bounded degree) and, in fact, it is
the main point of this paper to develop techniques for producing
analytic disks $C_{\alpha}$ in families. 

\smallskip\noindent{\slsf 6.} Corollary \ref{levi-gen} is proved in
section \ref{GEN}, where also a general position assumption is discussed.
Examples \ref{example1} and \ref{example2} are treated in section \ref{EXMP}.

\smallskip At the end I would like to give my thanks to the
Referee of this paper for the valuable remarks and suggestions.

\newsect[PINCH]{Extension to Pinched Domains}

\newprg[PINCH.pinch]{Analytic families and pinched domains}

By an {\slsf analytic family} of holomorphic mappings from $\Delta$
to $\Delta$ we understand the quadruple $(\calx,\pi,\cala,\Phi)$
where:

\begin{itemize}

\item $\calx$ is a complex manifold, which is either a finite
dimensional or a Banach one;

\smallskip\item a holomorphic submersion $\pi :\calx\to\cala$,
where  $\cala$ is a positive dimensional complex (Banach) manifold
such that for every $\alpha\in \cala$ the preimage $\calx_{\alpha}
:=\pi^{-1}(\alpha )$ is a disk;

\smallskip\item
a holomorphic map $\Phi :\calx\to \cc^2$ of generic rank $2$ such
that for every $\alpha\in \cala$ the image $\Phi (\calx_{\alpha})=
C_{\alpha}$ is a graph of a holomorphic function
$\phi_{\alpha}:\Delta \to \Delta$.
\end{itemize}
A family $(\calx , \pi, \cala, \Phi)$ we shall often call also an
analytic family of complex disks in $\Delta^2$. In our applications
$\cala$ will be always a neighborhood of some $\alpha_0$ and without
loss of generally we may assume for convenience that
$\phi_{\alpha_0}\equiv 0$, \ie that $C_{\alpha_0} = \Delta\times
\{0\}$. In this local case, after shrinking $\calx$ and $\Delta^2$
if necessary, we can suppose that $\calx = \Delta\times \cala$, and
we shall regard in this case $\Phi$ as a natural universal map

\begin{equation}
\eqqno(map-phi)
\Phi : (\lambda , \alpha) \to (\lambda , \phi_{\alpha}(\lambda))
\end{equation}
from $\Delta\times \cala$ to $\Delta^2$, writing $\Phi
(\lambda,\alpha) = (\lambda , \phi (\lambda ,\alpha))$ when
convenient, meaning $\phi (\lambda , \alpha) = \phi_{\alpha}
(\lambda)$. We shall often consider the case when $\cala$ is a
one-dimensional disk, in that case we say that our family is a
{\slsf complex one-parameter} analytic family. In this case taking
as $\cala$ a sufficiently small neighborhood of $\alpha_0$ and
perturbing $\d\Delta$ in $\lambda$-variable slightly we can suppose
without loss of generality that $\phi_{\alpha}$ doesn't vanish on
$\d\Delta$ if $\alpha\not=\alpha_0$. In particular the winding
number of $\phi_{\alpha}|_{\d\Delta}$ is constant for $\alpha\in
\cala\setminus \{\alpha_0\}$, see Proposition \ref{wind-fin-dim}.

\smallskip Denote as  $\bar\calp_{\calx,\Phi}$  the image  $\Phi (\calx)$,
where $(\calx, \pi , \Delta , \Phi)$ is some complex one-parameter
analytic family of complex disks in $\Delta^2$. Point $\lambda_0$
such that $\phi(\lambda_0,\alpha)\equiv 0$ as a function of $\alpha$
we call {\slsf a pinch} of $\bar\calp_{\calx ,\Phi}$ and say that
$\bar\calp_{\calx ,\Phi}$ has a pinch at $\lambda_0$.
Let us describe the shape of $\bar\calp_{\calx,\Phi}$ near a pinch
$\lambda_0$. Since $\phi (\lambda_0,\alpha)\equiv 0$ we can divide
it by $(\lambda - \lambda_0)^{l_0}$ with some (taken to be maximal)
$l_0\ge 1$. I.e., in a neighborhood of $(\lambda_0, \alpha_0)\in
\Delta\times \cala$ we can write
\begin{equation}
\eqqno(pinch1)
\phi (\lambda , \alpha) = (\lambda - \lambda_0)^{l_0}\phi_1(\lambda, \alpha),
\end{equation}
where $\phi_1(\lambda_0, \alpha)\not\equiv 0$. Set
\begin{equation}
\eqqno(map-phi1)
\Phi_1 : (\lambda , \alpha) \to (\lambda , \phi_1(\lambda , \alpha)).
\end{equation}
The image of $\Phi_1$ contains a bidisk $\Delta^2_r(\lambda_0,0)$ of some radius
$r>0$ centered at $(\lambda_0,0)$.
Therefore
\begin{equation}
\eqqno(pinch2)
\bar\calp_{\calx,\Phi} \supset \Delta^2_r(\lambda_0,0)\cap \{|z| < c |\lambda -
\lambda_0|^{l_0}\}
\end{equation}
with some constant $c>0$.

\begin{defi}
\label{pinch-d}
By a {\slsf pinched domain}  we shall understand an open neighborhood
$\calp$ of $\bar\Delta \setminus \Lambda$, where $\Lambda$ is a finite set of points in
$\Delta$, such that in a neighborhood of  every $\lambda_0\in \Lambda$ domain $\calp$
contains
\begin{equation}
\eqqno(pinch6)
 \Delta^2_r(\lambda_0,0)\cap \{|z| < c |\lambda - \lambda_0|^{l_0}\}\setminus
 \{(\lambda_0,0)\}.
\end{equation}
We shall call $l_0$ the {\slsf order of the pinch} $\lambda_0$.
\end{defi}
After shrinking $\Delta$ (in $\lambda$-variable) if necessary, we can
suppose that the set $\bar\calp_{\calx,\Phi}$ which corresponds to a complex
one-parameter analytic family $(\calx , \pi , \cala, \Phi)$ has only
finite number of pinches, say at $\lambda_1,...,\lambda_N$ of orders
$l_1,...,l_N$ respectively, and therefore $\calp_{\calx,\Phi}\deff\bar\calp_{\calx,\Phi}
\setminus \{\lambda_1,...,\lambda_N\}$ is a pinched domain. Remark that
$\bar\calp_{\calx,\Phi}$ obviously
contains every curve in a neighborhood $\calb$ of $\phi_0\equiv 0$ of the
subspace
\begin{equation}
\eqqno(pinch3)
\{ \phi \in \hol (\Delta , \Delta) : \ord_0(\phi , \lambda_j)\ge l_j\}
\subset \hol (\Delta , \Delta),
\end{equation}
which is of finite codimension.

\begin{rema} \rm
\label{hilb-str} {\bf (a)} Therefore, let us make the following
precisions: our pinched domains will be always supposed to have only
finitely many pinches and moreover, these pinches {\slsf do not
belong} to the corresponding pinched domain {\slsf by definition}.

\smallskip\noindent {\bf (b)} Hilbert manifold structure on $\calb$ (if needed)
can be insured by considering instead of $\hol (\Delta , \Delta)$
the Hilbert space $H^{1,2}_+(\sph^1)$ of Sobolev functions on the
circle, which holomorphically extend to $\Delta$, for example. This
will be done later in section \ref{N-TEST}. At that point it will be
sufficient for us to remark that extension along one-parameter
analytic families is equivalent to that of along of infinite
dimensional ones, and both imply the extension to pinched domains.
More precisely, the following is true:
\end{rema}

\begin{prop}
\label{pinch0} Let $(\calx , \pi , \cala ,\Phi)$ be a complex
one-parameter analytic family of complex disks in $\Delta^2$ and let
$\calp_{\calx,\Phi}$ be the corresponding pinched domain. Suppose
that a holomorphic function $f$ on $R_{1+\eps}$  meromorphically
extends along every $C_{\alpha}, \alpha\in\cala$. Let $\calb$ be the
{\slsf infinite dimensional} analytic family of complex disks in
$\Delta^2$ constructed as in \eqqref(pinch3). Then:

\smallskip\sli Function $f$ meromorphically extends to $\calp_{\calx , \Phi}$
as a function of two variables.

\smallskip\slii For every $\beta\in\calb$ the restriction $f|_{C_{\beta}
\cap R_{1+\eps}}$ extends to a meromorphic function on $C_{\beta}$
and the number of poles of these extensions is uniformly bounded.
\end{prop}
\proof Writing $\calx = \Delta\times \cala$ with $\cala = \Delta$ and
$\alpha_0=0$, and taking the preimage
$W\deff \Phi^{-1}(R_{1+\eps})$ in $\Delta\times \cala \equiv
\Delta^2$ we find ourselves in the following situation:

\smallskip\sli $W$ contains a ring domain (denote it by $W$ as well),
and $g\deff f\circ\Phi$ is meromorphic (or holomorphic, after
shrinking) on $W$.

\smallskip\slii For every $\alpha \in \cala $  the restriction $g|_{(\Delta\times\{\alpha\})\cap W}$
meromorphically extends to $\Delta\times \{\alpha\}$.

\smallskip The classical theorem of Levi, \cite{Lv, Si1}, implies now that $g$ meromorphically extends to
$\calx = \Delta\times \cala$, and this gives us the extension of $f$ to $\calp_{\calx , \Phi}$.

\smallskip For the proof of the extendability of $f|_{C_{\beta}\cap R_{1+\eps}}$ to $C_{\beta}$ for every
$\beta \in \calb$ close enough to zero let us first of all remark that $\Phi^{-1}(C_{\beta})$ is contained in a
relatively compact part of $\calx$. Indeed, take a pinch $\lambda_0$ and suppose without loss of generality
that $\lambda_0 = 0$. Write
\begin{equation}
\eqqno(weier0)
\phi (\lambda , \alpha) = \lambda^{l_0}\phi_1(\lambda
, \alpha)
\end{equation}
as in \eqqref(pinch1). Since $\phi_1(0,\alpha)\not\equiv 0$ we can use the Weierstrass preparation theorem
and present
\begin{equation}
\eqqno(weier1)
\phi_1(\lambda , \alpha) = u\cdot \left(\alpha^k + g_1(\lambda)\alpha^{k-1}+...+g_k(\lambda)\right)
\end{equation}
with $u(0,0)\not=0$ and $g_1(0)=...=g_k(0)=0$. Take the corresponding $\phi_{\beta}\in \calb$ with
graph $C_{\beta}$ and write it as $\phi_{\beta} (\lambda ) = c_0\lambda^{l_0}
\tilde\phi (\lambda)$. Consider the equation $\phi (\lambda ,\alpha) = \phi_{\beta} (\lambda )$, \ie
\begin{equation}
\eqqno(weier2)
\lambda^{l_0}\cdot \left(\alpha^k + g_1(\lambda)\alpha^{k-1}+...+g_k(\lambda)\right) = u^{-1} c_0\lambda^{l_0}
\tilde\phi (\lambda),
\end{equation}
or, equivalently
\begin{equation}
\eqqno(weier3)
\alpha^k + g_1(\lambda)\alpha^{k-1}+...+g_k(\lambda) = u^{-1} c_0\tilde\phi (\lambda).
\end{equation}
For $\lambda \sim 0$ all solutions $\alpha_1(\lambda),...,\alpha_k(\lambda)$ of \eqqref(weier3) are close to zero,
provided $c_0$ is small enough. This proves our assertion that $\Phi^{-1}(C_{\beta})\comp \calx$ and implies
that $f|_{C_{\beta}\cap R_{1+\eps}}$ meromorphically extends to $C_{\beta}$.

\smallskip The orders of poles of meromorphic function $g$ of two variables $(\lambda , \alpha)$ is bounded
on every relatively compact part of $\calx = \Delta\times \cala $ and therefore the orders of
poles of our extensions are also bounded.

\smallskip\qed

\begin{rema} \rm
\label{inf-pinch}
Remark that $f$ meromorphically extends to the pinched domain $\calp_{\calb}$ swept
by the family $\calb$ as well, simply because it is the same domain (up to shrinking).
\end{rema}

\newprg[PINCH.proof]{Proof of Theorem \ref{levi-nonlin}}
We start with the proof of item (\slii  first. Without loss of
generality we may assume that $\phi_0\equiv 0$.  Indeed, the 
condition $\phi_0\equiv 0$ is not a restriction
neither here, nor anywhere else in this paper, because it can be
always achieved by the coordinate change
\begin{equation}
\eqqno(coord-ch)
\begin{cases}
\lambda \to \lambda ,\cr
z\to z-\phi_{0}(\lambda).
\end{cases}
\end{equation}
Furthermore, when considering the extension of a meromorphic function 
$f$ from a ring
domain $R_{1+\eps}$ to the bidisk $\Delta_{1+\eps}\times \Delta$ one
can suppose that $f$ is holomorphic on $R_{1+\eps}$ (after shrinking
$R_{1+\eps}$ if necessary and after multiplying by some power of $z$),
and moreover, decomposing $f = f^+ + f^-$ where $f^+$ is holomorphic
in $\Delta_{1+\eps}\times \Delta$ and $f^-$ in $(\pp^1\setminus
\bar\Delta)\times \Delta$, one can subtract $f^+$ from $f$ and
suppose that $f^+\equiv 0$. That means that we can suppose that $f$
has the Taylor decomposition
\begin{equation}
\eqqno(taylor0)
f(\lambda , z) = \sum_{n=0}^{\infty} A_n(\lambda) z^n
\end{equation}
in $R_{1+\eps}$ with
\begin{equation}
\eqqno(taylor-A)
A_n(\lambda) = \sum_{l=-\infty}^{-1}a_{n,l}\lambda^l.
\end{equation} 
As the result along this proof we may suppose that $f=f^-$ and $f^-$
is holomorphic in $A_{1-\eps, 1+\eps}\times\Delta_{1+2\eps}$.
Therefore for $|\lambda|$ near $1$ the Taylor expansion  of $f$ writes
as
\begin{equation}
\eqqno(taylor)
f(\lambda ,z)= \sum_{n=0}^{\infty}\frac{1}{n!}\frac{\d^n f(\lambda , 0)}{\d z^n}z^n =
\sum_{n=0}^{\infty}A_n(\lambda)z^n,
\end{equation}
and we have the estimates
\begin{equation}
\eqqno(cauchy0)
\left| A_n(\lambda)\right| = \frac{1}{n!}\left\vert \frac{\d^n f(\lambda ,0)}{\d z^n}\right\vert \le
\frac{C}{(1+\eps )^n},
\end{equation}
for some constant $C$, all $k\in \nn$ and all $\lambda\in\sph^1\deff \d\Delta$.
Under the assumptions of the Theorem we see that meromorphic extensions $f_k(\lambda)$ of
$f(\lambda , \phi_k(\lambda))$ have uniformly bounded number of poles counted with multiplicities.
As well as the numbers of zeroes of $\phi_k$ are uniformly bounded too. Up to taking a subsequence
we can suppose that:

\smallskip{\slsf a)} The number of poles of $f_k$-s, counted with multiplicities, is constant, say $M$,
and these poles converge to the finite set $b_1,...,b_M\in \Delta_{1-\eps}$ with corresponding
multiplicities, \ie some of $b_1,...,b_M$ may coincide.

\smallskip{\slsf b)} The number of zeroes of $\phi_k$, counted with multiplicities, is also constant, 
say $N$ and these zeroes converge to a finite set with corresponding multiplicities. We shall denote 
it as $a_1,...,a_N$, meaning that some of them can coincide.

\smallskip\noindent{\slsf Step 1.} For every $k$ take a Blaschke product $P_k$ having zeroes exactly 
at poles of $f_k$ with corresponding multiplicities and subtract from $\{P_k\}$
a converging subsequence with the limit
\begin{equation}
\eqqno(blasch1)
P_0(\lambda) = \prod_{i=1}^M\frac{\lambda - b_i}{1-\bar b_i\lambda}.
\end{equation}
Holomorphic functions $g_k\deff P_kf_k$ have uniformly bounded 
modulus  on $\Delta$ and converge to some $g_0$, 
with modulus bounded by $C$ (a constant from \eqqref(cauchy0)). Therefore
$f_k$ converge on compacts of $\Delta\setminus \{b_1,...,b_M\}$ to a
meromorphic function, which is nothing but $A_0$, and it satisfies
the estimate
\begin{equation}
\eqqno(est-0) |A_0(\lambda)|\le \frac{CC_1}{|\lambda -
b_1|...|\lambda -b_M|},
\end{equation}
where $C_1 = \max\{\Pi_{i=1}^M|1-\bar b_i\lambda| :|\lambda|\le
1\}$.

\smallskip\noindent{\slsf Step 2.} Consider the function
\begin{equation}
\eqqno(funct-1) f_1(\lambda ,z) \deff \frac{f(\lambda , z) -
A_0(\lambda )}{z},
\end{equation}
and the following functions
\begin{equation}
\eqqno(rasnost1) f_{1,k}(\lambda) \deff f_1(\lambda
,\phi_k(\lambda)) = \frac{f(\lambda, \phi_k(\lambda)) - A_0(\lambda
)}{\phi_k(\lambda)}.
\end{equation}
These functions are well defined and meromorphic on
$\Delta_{1+\eps}$, the equality in \eqqref(rasnost1) has sense on
$A_{1-\eps,1+\eps}$. After taking a subsequence we see that poles of
$f_{1,k}$, which are different from zeroes of $\phi_k$, converge to the
same points $b_1,...,b_N$. So the
multiplicities do not increase. Let again $P_k$ be the Blaschke
product having zeroes at poles of $f_{1,k}$ with corresponding multiplicities.
After taking a subsequence $P_k$ uniformly converge to a
corresponding Blaschke product $P_0$ and holomorphic functions
$g_k\deff P_kf_{1,k}$ uniformly converge to some holomorphic
function $g_0$. In $A_{1-\eps , 1+\eps}$ it is straightforward to see that
$g_0 = P_0A_1$. This proves that $A_1$ (if not identically zero),
has at most $N+M$ poles counting with multiplicities and these poles
are located at $a_1,...,a_N,b_1,...,b_M$.

Moreover, for $|\lambda|=1$ from \eqqref(cauchy0) we have the estimate
\begin{equation}
\eqqno(cauchy1)
|P_0(\lambda)A_1(\lambda)| \le \frac{C}{1+\eps},
\end{equation}
which implies the estimate
\begin{equation}
\eqqno(est-1) |A_1(\lambda)| \le \frac{1}{|\lambda -a_1|...|\lambda
-a_N||\lambda -b_1||\lambda - b_N|}\cdot \frac{CC_1C_2}{1+\eps}
\end{equation}
for $\lambda\in \Delta\setminus \{a_1,...,a_N,b_1,...,b_M\}$. Here
$C_1 = \max\{\Pi_{i=1}^N|1-\bar a_i\lambda| :|\lambda|\le 1\}$.
Denote from now $CC_1C_2$ by $C'$.

\smallskip\noindent{\slsf Step 3.} Suppose we proved that $A_n$ extends to a meromorphic function in
$\Delta$ with the estimate
\begin{equation}
\eqqno(est-n) |A_n(\lambda)| \le \frac{1}{\prod_{j=1}^N|\lambda
-a_j|^n\prod_{j=1}^{M}|\lambda - b_j|}\cdot \frac{C'}{(1+\eps )^n}
\end{equation}
for $\lambda\in \Delta\setminus \{a_1,...,a_N,b_1,...,b_M\}$. Remark that \eqqref(est-n) means,
in particular, that $A_0,...,A_n$ have no other poles than $a_1,...,b_N$ with corresponding multiplicities.
Apply considerations as above to
\[
f_{n+1}(\lambda , z) = \frac{1}{z^{n+1}}\left(f(\lambda , z) -
\sum_{j=0}^nA_j(\lambda)z^j\right),
\]
\ie  consider
\[
f_{n+1, k}(\lambda) = \frac{1}{\phi_k^{n+1}}\left(f(\lambda ,
\phi_k) - \sum_{j=0}^nA_j(\lambda)\phi_k^j\right)
\]
and repeat the same demarche with Blaschke products. Remark only that products $A_j(\lambda)\phi_k^j$ have no
poles at zeroes of $\phi_k$. On the boundary $\{|\lambda |=1\}$ functions $|f_{n+1, k}(\lambda)|$ are bounded
by $C/(1+\eps)^{n+1}$
due to Cauchy inequalities and therefore we get the conclusion that $A_{n+1}$ meromorphically extends to $\Delta$
with the estimate
\begin{equation}
\eqqno(est-n+1) |A_{n+1}(\lambda)| \le
\frac{1}{\prod_{j=1}^N|\lambda -a_j|^{n+1}\prod_{j=1}^M|\lambda -
b_j|}\cdot \frac{C'}{(1+\eps )^{n+1}}.
\end{equation}
Estimate \eqqref(est-n+1) implies that \eqqref(taylor) converges in the domain
\begin{equation}
\eqqno(pinch4)
\left\{(\lambda ,z)\in \Delta^2: |z|< c |\lambda - a_{j_1}|^{l_1}...|\lambda - a_{N_1}|^{l_{N_1}}
\right\}\setminus \bigcup_{i=1}^{M}\{\lambda =b_{i}\},
\end{equation}
for an appropriately chosen $c>0$. Here $N_1$ is the number of
{\slsf different} $a_j$-s, which are denoted as
$a_{j_1},...,a_{N_1}$ having corresponding multiplicities
$l_{1},...,l_{N_1}$. In particular we mean here
that $b_{i}$ are different from $a_{j_1}$ for all $i,j_1$. Estimate \eqqref(est-n+1) implies that the
extension of $f\cdot \prod_{j=1}^M(\lambda - b_j)$ to
\eqqref(pinch4) is locally bounded near every vertical disk
$\{\lambda =b_{i_1}\}$ and therefore extends across it by Riemann
extension theorem. We conclude that $f$ extends as a meromorphic
function to the pinched domain
\begin{equation}
\eqqno(pinch5)
\calp  = \left\{(\lambda ,z)\in \Delta^2: |z|< c |\lambda - a_{j_1}|^{l_1}...|\lambda - a_{N_1}|^{l_{N_1}}
\right\} ,
\end{equation}
and this proves the part (\slii of Theorem \ref{levi-nonlin}.

\smallskip\noindent (\sli  Take now any holomorphic function $\phi$ in
$\Delta_{1+\eps}$ of the form
\[
\phi (\lambda ) = (\lambda -a_{j_1})^{l_1}...(\lambda - a_{N_1})^{l_{N_1}}\psi
\]
with $\psi$ small enough in order that the graph $C_{\phi}$ is
contained in $\calp$ (more precisely should be $C_{\phi}\cap
\left(\Delta\setminus \{a_{j_1},...,a_{N_1}\}\right)\times \Delta
\subset \calp $). To prove the part (\sli of our theorem we need to
prove the following

\smallskip\noindent{\slsf Step 4.} {\it $f(\lambda , \phi(\lambda))$
meromorphically extends from $A_{1-\eps, 1+\eps}$ to
$\Delta_{1+\eps}$.} Indeed, $f(\lambda , \phi (\lambda))$ is
meromorphic on $\Delta\setminus \{a_{j_1},...,a_{N_1}\}$. At the
same from the estimate \eqqref(est-n) we see that the terms in the
series
\begin{equation}
\eqqno(sweep) f(\lambda , \phi (\lambda)) = \sum_{n=0}^{\infty}
A_n(\lambda)\phi ^n(\lambda)
\end{equation}
are, in fact, holomorphic in a neighborhood of every $a_j$ and
converge normally there, provided $\norm{\psi}_{\infty}$ was taken
small enough. Uniform boundedness of the number of poles follows now from
Proposition \ref{pinch0}. Part (\sli is proved.

\smallskip\qed

\begin{rema} \rm
\label{pr-levi-c}
In order to prove Corollary \ref{levi-c} remark
that pinches that appeared along the proof of Theorem \ref{levi-nonlin} are
limits of zeroes of $\phi_k$. General position assumption means
that for every $\lambda_0\in \Delta$ we can take a subsequence such
that the resulting pinched domain will not have a pinch in
$\lambda_0$. The rest follows.
\end{rema}

\newsect[N-TEST]{Extension Along Finite Dimensional Families}

\newprg[N-TEST.sing]{Properties of the Singular Integral Transform}

By $L^{1,2}(\sph^1)$ we denote the Sobolev space of complex valued
functions on the unit circle having their first derivative in $L^2$. This
is a complex Hilbert space with the scalar product $(h,g) = \int_{0}^{2\pi }
[h(e^{i\theta})\bar g(e^{i\theta}) + h'(e^{i\theta})\bar g'(e^{i\theta})]d\theta$.
Recall that by Sobolev Imbedding Theorem $L^{1,2}(\sph^1)\subset
\calc^{\frac{1}{2}}(\sph^1)$, where $\calc^{\frac{1}{2}}(\sph^1)$ is
the space of H\"older $\frac{1}{2}$-continuous functions on
$\sph^1$.

\smallskip For the convenience of the reader we recall few well known facts
about the Hilbert Transform in $L^{1,2}(\sph^1)$.

\begin{lem}
\label{oper-p-lem}
A function $\phi\in L^{1,2}(\sph^1)$ extends holomorphically to $\Delta$
if and only if the following condition is satisfied:
\begin{equation}
\eqqno(oper-p)
P(\phi)(\tau) := \frac{-1}{2\pi
i}\int_{\sph^1}\frac{\phi(t)-\phi(\tau)}{t-\tau}dt \equiv 0.
\end{equation}
\end{lem}

\proof The fact that $\phi$ extends holomorphically to $\Delta$ can be obviously expressed as

\begin{equation}
\eqqno(oper-p1)
\lim_{z\to \tau ,z\in \Delta}\frac{1}{2\pi i}\int_{\sph^1}\frac{\phi(t)}{t-z}dt = \phi(\tau )
\end{equation}
for all $\tau\in\sph^1$.  Write then

\begin{equation}
\eqqno(oper-p2)
\lim_{z\to \tau ,z\in \Delta}\frac{1}{2\pi
i}\int_{\sph^1}\frac{\phi(t)}{t-z}dt = \lim_{z\to \tau ,z\in \Delta} \frac{1}{2\pi
i}\int_{\sph^1}\frac{\phi(t)-\phi(\tau )}{t-z}dt +
\end{equation}
\[
+ \lim_{z\to \tau
,z\in \Delta}\frac{1}{2\pi i}\int_{\sph^1}\frac{\phi(\tau)}{t-z}dt
= - P(\phi)(\tau ) + \phi(\tau ).
\]
From \eqqref(oper-p1) and \eqqref(oper-p2) we immediately get \eqqref(oper-p).

\smallskip\qed

\smallskip
Denote by $\sph^1_{\eps}(\tau)$ the circle $\sph^1$ without the
$\eps$-neighborhood of $\tau$. Consider the following singular
integral operator (the Hilbert Transform)
\begin{equation}
\label{oper-s}
S(\phi)(\tau) :=  {\mathbf{p.v.}} \frac{1}{\pi
i}\int_{\sph^1}\frac{\phi(t)}{t-\tau}dt := \lim_{\eps\to
0}\frac{1}{\pi i} \int_{\sph^1_{\eps}(\tau)}\frac{\phi(t)}{t-\tau}dt .
\end{equation}
In the sequel we shall write simply
\begin{equation}
\label{oper-ss}
S(\phi)(\tau) :=  \frac{1}{\pi
i}\int_{\sph^1}\frac{\phi(t)}{t-\tau}dt ,
\end{equation}
\ie the integral in the right hand side will be always understood
in the sense of the principal value.

\begin{lem} The following relation between operators $S$ and $P$ holds
\begin{equation}
\eqqno(id-p-s)
S = -2P + \id .
\end{equation}
Therefore a function $\phi\in L^{1,2}(\sph^1)$ holomorphically
extends to the unit disk if an only if
\begin{equation}
\eqqno(hilb)
S(\phi)(\tau) \equiv \phi(\tau).
\end{equation}
\end{lem}
\proof Write

\[
\frac{1}{2}S(\phi)(\tau) = \frac{1}{2\pi i}\int_{\sph^1}\frac{\phi(t)}{t-\tau}dt = \lim_{\eps\to
0}\frac{1}{2\pi i} \int_{\sph^1_{\eps}(\tau)}\frac{\phi(t)}{t-\tau}dt
=
\]
\[
= \frac{1}{2\pi i}\int_{\sph^1}\frac{\phi(t)-\phi(\tau )}{t-\tau} +
\lim_{\eps\to 0}\frac{1}{2\pi i}\int_{\sph^1_{\eps}(\tau)}\frac{\phi(\tau)}{t-
 \tau}dt = - P(\phi)(\tau) + \frac{1}{2}\phi(\tau).
\]
Therefore one has
\begin{equation}
S(\phi) = -2 P(\phi) + \phi ,
\end{equation}
which is \eqqref(id-p-s), and which implies \eqqref(hilb).

\smallskip\qed

\smallskip
Denote by $H^{1,2}_+(\sph^1)$  the subspace of $L^{1,2}(\sph^1)$
which consists of functions holomorphically extendable to the unit
disk $\Delta$. By $H^{1,2}_-(\sph^1)$ denote the subspace of
functions holomorphically extendable to the complement of the unit
disk in the Riemann sphere $\pp^1$ and zero at infinity. Observe
the following orthogonal decomposition
\begin{equation}
\eqqno(ort-decomp)
L^{1,2}(\sph^1) = H^{1,2}_+(\sph^1) \oplus H^{1,2}_-(\sph^1).
\end{equation}
We finish this review with the following:
\begin{lem}
\label{proj-p}
\sli $P$ and $S$ are bounded linear operators on $L^{1,2}(\sph^1)$ and
\begin{equation}
\eqqno(oper-s2)
S^2=\id.
\end{equation}
\slii Moreover, on the space $H^{1,2}_+(\sph^1)$ operator $S$ acts as identity
and on the space $H^{1,2}_-(\sph^1)$ as $-\id$.

\smallskip\noindent\sliii Consequently $P$ is an orthogonal projector
onto $H^{1,2}_-(\sph^1)$.
\end{lem}
For the proof of \eqqref(oper-s2) we refer to \cite{MP} pp. 46, 50,
69. In fact, since $S = -2P + \id $ and because of $\ker P =
H^{1,2}_+(\sph^1)$, we see that $S=\id $ on $H^{1,2}_+(\sph^1)$.
From \eqqref(oper-s2) and \eqqref(id-p-s) we also see that $P= \id $
on $H^{1,2}_-(\sph^1)$, \ie $P$ projects $L^{1,2}(\sph^1)$ onto
$H^{1,2}_-(\sph^1)$ parallel to $H^{1,2}_+(\sph^1)$.

\smallskip This lemma clearly implies the following
\begin{corol}
\label{mer-ext}
Function $\phi\in L^{1,2}(\sph^1)$ extends to a meromorphic function in
$\Delta$ with not more than $N$ poles if and only if $P(\phi)$ is rational, 
zero at infinity and has not more than $N$ poles.
\end{corol}
Indeed, decompose $\phi = \phi^+ + \phi^-$ according to \eqqref(ort-decomp).
$\phi$ is meromorphic with at most $N$ poles, all in $\Delta$, if and only 
if $\phi^-$ is such. Which means that $\phi^-$ should be ratioanl with at 
most $N$ poles. But since, according to (\sliii of Lemma \ref{proj-p} one has
$P(\phi) = \phi^-$, the last is equivalent to the fact that $P(\phi)$ 
is rational with at most $N$ poles.

\newprg[N-TEST.fin-dim]{Case of finite dimensional families}

To clarify the finite vs. infinite dimensional issues in this
paper let us give a simple proof of Theorem \ref{levi-nonlin}
in the special case when $\phi_k$ belong to
an analytic family $\{\phi_{\alpha}\}_{\alpha\in\cala}$
parameterized by a finite dimensional complex manifold $\cala$. More
precisely, as in section \ref{PINCH}, we are given a complex
manifold $\calx$, a holomorphic submersion $\pi :\calx\to\cala$ such
that for every $\alpha\in \cala$ the preimage $\calx_{\alpha}
:=\pi^{-1}(\alpha)$ is a disk. We are given also a holomorphic map
$\Phi:\calx\to \cc^2$ such that for every $\alpha\in \cala$ the
image $\Phi(\calx_{\alpha})= C_{\alpha}$ is a graph of a holomorphic
function $\phi_{\alpha}:\Delta_{1+\eps}\to \Delta $. We
shall regard $\cala$ as a (locally closed) complex submanifold of
$H^{1,2}_+(\sph^1)$. And, finally, by saying that
$\{\phi_k\}_{k=1}^{\infty}$ belong to
$\{\phi_{\alpha}\}_{\alpha\in\cala}$ we mean that there exist
$\alpha_k\in \cala , \alpha_k\to \alpha_0\in\cala$, such that
$\phi_k = \phi_{\alpha_k}$ for $k\ge 0$.

\smallskip After shrinking, if necessary, we can suppose that our
function $f$ is holomorphic on $R_{1+\eps}=A_{1-\eps,1+\eps}\times
\Delta_{1+\eps}$. Consider the following analytic mapping
$F:L^{1,2}(\sph^1)\to L^{1,2}(\sph^1)$
\begin{equation}
\eqqno(map-f)
F : \phi (\lambda) \to f(\lambda , \phi (\lambda)),
\end{equation}
and consider also the following integral operator $\calf
:H^{1,2}_+(\sph^1)\to H_-^{1,2}(\sph^1)$
\begin{equation}
\eqqno(map-F) \calf (\phi)(\lambda) = \frac{-1}{2\pi
i}\int_{\sph^1}\frac{f(\zeta,\phi (\zeta)) - f(\lambda,\phi
(\lambda))}{\zeta - \lambda}d\zeta = P\left(F(\phi )\right).
\end{equation}
According to Lemma \ref{oper-p-lem} $f(\lambda,\phi (\lambda))$
extends to a holomorphic function in $\Delta_{1+\eps}$ if and only
if $\calf (\phi)=0$, and according to Corollary \ref{mer-ext}
it extends meromorphically to $\Delta_{1+\eps}$ with at most $N$
poles in $\Delta_{1-\eps}$ if an only if $\calf  (\phi)$ is a 
boundary value of a rational function with at most $N$ poles all 
in $\Delta_{1-\eps}$.

\begin{thm}
\label{fin-dim-case}
Let $f$ be a meromorphic function on $R_{1+\eps}$ and
$\{\phi_k:\Delta_{1+\eps}\to \Delta\}_{k=1}^{\infty}$ be a sequence
of holomorphic functions converging to some $\phi_0 : \Delta_{1+\eps}
\to \Delta$, $\phi_k\not\equiv \phi_0$ for all $k$. Suppose that:

\smallskip a) $\{\phi_k\}$ belong to a finite dimensional analytic 
family $\{\phi_{\alpha}\}_{\alpha\in \cala}$, \ie $\phi_k=\phi_{\alpha_k}$
for some 

\quad $\alpha_k\in \cala$ and $\alpha_k\to \alpha_0$ in $\cala$ with 
$\phi_0=\phi_{\alpha_0}$;

\smallskip b) for every $k$ the restriction
$f|_{C_{k}\cap R_{1+\eps}}$ is well defined and extends to a meromorphic

\quad function on the curve $C_{k}$;

\smallskip c) the number of poles counting with
multiplicities of these extensions is uniformly

\quad bounded.

Then there exists a complex disk $\Delta\subset\cala$ containing
$\alpha_0$ such that for every $\alpha\in \Delta$ the restriction
$f|_{C_{\alpha}\cap R_{1+\eps}}$ meromorphically extends to
$C_{\alpha}$, and the number of poles of these extensions counting
with multiplicities is uniformly bounded.
\end{thm}
\proof \sli Consider the holomorphic case first. Restrict $\calf$ to
$\cala$ to obtain a holomorphic map $\calf_{\cala}:\cala\to
H^{1,2}_-(\sph^1)$. $f|_{C_{\alpha}\cap R_{1+\eps}}$ holomorphically
extends to $C_{\alpha}$ if and only if $\calf (\alpha)=0$. Therefore
we are interested in the zero set $\cala^0$ of $\calf_{\cala}$. But
the zero set $\cala^0$ of a holomorphic mapping from a finite
dimensional manifold is a {\slsf finite dimensional} analytic set.
Since this set contains a converging sequence $\{\alpha_k\}$ it has
positive dimension. Denote (with the same letter) by $\cala^0$ a
positive dimensional irreducible component of our zero set which
contains an infinite number of $\phi_{\alpha_k}$-s. Suppose, up to 
replacing $\alpha_k$ by a subsequence, that all $\alpha_k$ are in 
$\cala^0$. Let $\calx^0$ be the corresponding universal family, \ie 
the restriction of $\pi:\calx\to \cala$ to $\cala^0$, and
$\Phi^0:\calx^0 \to \Delta_{1+\eps}\times\Delta$ the corresponding
evaluation map. $\Phi^0$ should be of generic rank two, otherwise
$\phi_k$ would be constant. Therefore $\cala^0$ contains a complex disk
through $\alpha_0$ with properties as required.

\smallskip\noindent\slii
The meromorphic extension in this case is also quite simple. Without
loss of generality we suppose that all extensions $f_{\alpha_k}(\lambda)$
have at most $N$ poles counting with multiplicities. Since
$f_{\alpha_k}(\lambda)=f(\lambda ,\phi_{\alpha_k}(\lambda ))$ for
$\lambda\in A_{1-\eps,1+\eps}$, all poles of these extensions are
contained in $\bar\Delta_{1-\eps}$. Denote by $R^N(1-\eps)$ the
subset of $H^{1,2}_-(\sph^1)$ which consists of
rational functions, holomorphic on $\pp^1\setminus \bar\Delta$,
zero at infinity and having not more than $N$ poles, all contained in
$\bar\Delta_{1-\eps}$. $R^N(1-\eps)$ can be explicitly described as
the set of the following functions:
\begin{equation}
\eqqno(RN-eps)
R^N(1-\eps) = \left\lbrace \sum_{j}(z-a_j)^{-m_j}\sum_{k=0}^{m_j-1}c_{jk}(z-a_j)^k:
c_{jk}\in \cc ,
a_j\in \bar\Delta_{1-\eps}, \sum_jm_j=N\right\rbrace .
\end{equation}
Let us note that $\calf_{\cala} (\phi_{\alpha_k})\in R^N(1-\eps)$ for all
$k$ and that the set $\cala^N$ of those $\alpha\in\cala$ that
$f(\lambda ,\phi_{\alpha}(\lambda))$ is meromorphically extendable
to $\Delta$ with not more $N$ poles, all in $\bar\Delta_{1-\eps}$, is in fact
$\calf_{\cala}^{-1}(R^N(1-\eps ))$.

\smallskip Set $g_0=\calf_{\cala} (\phi_{\alpha_0})$. From \eqqref(RN-eps)
we see that $R^N(1-\eps)$ is a finite dimensional subspace of $
H^{1,2}_-(\sph^1)$. Therefore we can take an orthogonal complement 
$H\subset H^{1,2}_-(\sph^1)$ to it at $g_0$ in such a way that
$H^{1,2}_-(\sph^1) = R^N(1-\eps)\times H$ locally in a neighborhood
of $g_0$. Denote by $\Psi$ the composition of $\calf_{\cala} $ with  the
projection onto $H$. Now $\cala^N$ is the zero set of $\Psi $ and
therefore we are done as in the case (\sli .

\medskip\qed

Let us make a few remarks concerning the finite dimensional case of the
last theorem.

\begin{rema} \rm
\label{fin-dim} {\bf a)} Horizontal disks belong to an one
dimensional family, non-horizontal straight disk to two-dimensional.
Therefore Theorem \ref{fin-dim-case} generalizes Hartogs-Levi
theorem and the result of Dinh.

\smallskip\noindent{\bf b)} We do not claim in Theorem \ref{fin-dim-case} that the disks $\Delta$
contains $\alpha_k$ for $k>>1$ and this is certainly not true in general. What is true
is that the set $\cala^0$ (or $\cala^N$) of {\slsf all} $\alpha \in \cala$ such that
$f$ extends along $C_{\alpha}$ is an analytic set of positive dimension, so contains an
analytic disk with center at $\alpha_0$.
\end{rema}

\begin{rema} \rm
At the same time Theorem \ref{fin-dim-case} is a particular case of Theorem
\ref{levi-nonlin} because of the following observation.

\begin{prop}
\label{wind-fin-dim}
Let $(\calx , \pi, \cala, \Phi)$ be a finite dimensional analytic family of holomorphic
maps $\Delta_{1+\eps}\to \Delta$ and let $\alpha_0\in \cala$ be a point. Then
there exists a neighborhood $V\ni \alpha_0$, a complex hypersurface $A$ in $V$
and a radius $r\sim 1$ such that for $\alpha \in V\setminus A $ the restriction
$(\phi_{\alpha} - \phi_0)|_{\d\Delta_r}$ doesn't vanish and therefore
$\var\arg_{\d\Delta}(\phi_{\alpha} - \phi_{\alpha_0})$ is constant on $V\setminus A$.
If, in particular, $(\calx , \pi, \cala, \Phi)$ is a one-parameter family then
$A=\{\alpha_0\}$.
\end{prop}
\proof After shrinking we can suppose that $\calx = \Delta_{1+\eps}\times \Delta^n$,
$\alpha_0 = 0$ and $\phi_0 \equiv 0$. Since $\Phi : \calx \to \Delta_{1+\eps}\times
\Delta$ writes as $(\lambda , \alpha)\to (\lambda ,\phi (\lambda,\alpha))$ we can
consider the zero divisor $\calz = \phi^{-1}(0)$ of $\phi$. $\calz$ is not empty,
because $\phi (\lambda , 0) \equiv 0$, and is proper, because $\Phi$ is of generic
rank two. Denote by $\calz_1$ the union of all irreducible components of $\calz$
which  contain $\Delta_{1 + \eps}\times \{0\}$. Set $A\deff \calz_1\cap \{0\}$
and remark that $A$ is a hypersurface in $\Delta^n$.

\smallskip Denote by $\calz_0$ the union of all irreducible components of $\calz$
which do not contain $\Delta_{1+\eps}\times \{0\}$. Intersection $\calz_0\cap
\Delta_{1+\eps}\times \{0\}$ is a discrete set. Therefore we can find $r\sim 1$
such that $\calz_0\cap \d\Delta_r = \emptyset$. Now it is clear that for a
sufficiently small neighborhood $V\ni 0$ in $\Delta^n$ we have that $\phi (\cdot , \alpha)$
doesn't vanish on $\d\Delta_r$ provided $\alpha \in V\setminus A$. Then
$\var\arg_{\d\Delta}(\phi_{\alpha})$ is clearly constant. In one-parameter
case $A$ is discrete but contains $\alpha_0$.

\smallskip\qed

\smallskip\noindent{\bf c)} Let us remark that in general test sequence
doesn't belong to any finite dimensional family. Take for example
$\phi_k(\lambda) = \frac{1}{k}\lambda^2 + e^{-k}\lambda^k$. Therefore Theorem
\ref{levi-nonlin} properly contains Theorem \ref{fin-dim-case}.
\end{rema}

\begin{rema} \rm
\label{tot-real} {\bf a)} If $C_k$ are intersections of
$\Delta_{1+\eps}\times \Delta$ with algebraic curves of bounded degree, 
then they are included in a finite dimensional analytic (even algebraic
in this case) family.

\smallskip\noindent{\bf b)} If $\phi_k(\d\Delta)\subset M$, where
$M$ is totally real in $\d\Delta\times\bar\Delta$, and have bounded
Maslov index then they are included in a finite dimensional analytic
family.

\smallskip\noindent{\bf c)}
If we do not suppose {\slsf ad hoc} that $\phi_k$ belong to some
finite dimensional analytic family of holomorphic functions then the
argument above is clearly not sufficient. The following example is 
very instructive. Consider a holomorphic map $\calf :l^2\to l^2\oplus
l^2$ defined as
\begin{equation}
\eqqno(coeure)
\calf  : \{z_k\}_{k=1}^{\infty} \to \left\{\{z_k(z_k-1/k)\}\oplus
\{z_kz_j\}_{j>k}\right\}.
\end{equation}
The zero set of $\calf$ is a sequence $\{Z_k = (0,...,0,1/k,0,...)
\}_{k\ge 1}\subset l^2$ together with zero. These $Z_k$-s might well 
be ours $\phi_k$-s and therefore we cannot conclude the existence of
families in the zero set of our $\calf$ from \eqqref(map-F) at this 
stage.

\smallskip\noindent{\bf d)} Example \ref{example1} has precisely the 
feature as above with $\calf$ being the integral operator \eqqref(map-F).
\end{rema}

\newsect[GEN]{Mappings to Complex Spaces}

In this section we shall prove Corollary \ref{levi-gen}. The proof
consists in making a reduction to the  holomorphic function case of
Theorem \ref{levi-nonlin}. This reduction will follow the lines of
arguments developed in \cite{Iv1,Iv2, Iv3, Iv4}. For the convenience of
the reader we shall briefly recall the key statements from these
papers which are relevant to our present task.

\newprg[GEN.cont]{Continuous families of analytic disks}

Analytic disk in a complex space $X$ is a holomorphic map
$h:\Delta\to X$ continuous up to the boundary. Recall that a complex
space $X$ is called disk-convex if for every compact $K\comp X$
there exists another compact $\hat K$ such that for every analytic
disk $h:\bar\Delta \to X$ with $h(\d\Delta)\subset K$ one has
$h(\bar\Delta)\subset \hat K$. $\hat K$ is called the {\slsf disk
envelope} of $K$. All compact, Stein, $1$-convex complex spaces are
disk-convex. 

\smallskip Given a meromorphic mapping $f:R_{1+\eps}\to X$, where
$R_{1+\eps} = A_{1-\eps, 1+\eps}\times \Delta $, we
can suppose without loss of generality that $f$ is holomorphic on
$R_{1+\eps}$ and that $f(R_{1+\eps})$ is contained in some compact
$K$. We suppose that our space $X$ is reduced and that it
is equipped with some Hermitian metric form $\omega$. Denote by $\nu
=\nu (\hat K)$
the minima of areas of rational curves in the disk envelope
$\hat K$ of $K$. Remark that $\nu$ is achievable by some rational curve
and therefore $\nu >0$. We are given an {\slsf uncountable} family of 
disks $\{C_t: t\in T\}$ which are the graphs of holomorphic functions
$\phi_t:\Delta_{1+\eps}\to \Delta$. Remark that the condition on 
our family of disks to be {\slsf test} (see Introduction) implies, in 
particular, that they are all distinct. In the sequel when writing 
$f(C_t)$ we mean more precisely $f|_{C_t}(C_t)$, \ie the restriction 
$f$ to $C_t$.  This is an analytic disk in $X$ and since $f|_{C_t}(\d
C_t)\subset K$ we see that for every $t\in T$ one has $f(C_t)\subset 
\hat K$. For every natural $k$ set
\begin{equation}
\eqqno(area1)
T_k =\left\{t\in T: \area\left(\Gamma_{f|_{C_t}}\right)\le
k\frac{\nu}{2}\right\},
\end{equation}
where $\Gamma_{f|_{C_t}}=:\Gamma_t$ is the graph of $f|_{C_t}$ in
$\Delta^2_{1+\eps}\times X$ and area is taken with respect to the
standard Euclidean form $\omega_e=dd^c(|\lambda|^2+|z|^2)$ on
$\cc^2$ and $\omega$ on $X$. For some $k$ the set $T_k\setminus
T_{k-1}$ is uncountable, so denote this set as $T$ again.

\smallskip It will be convenient in the sequel to consider our parameter
space $T$ as a subset of the space of $1$-cycles in $\Delta^2_{1+\eps}\times X$.
Let us say few words about this issue. For general facts about cycle spaces
we refer to \cite{Ba}, for more details concerning our special situation
to \S 1 of \cite{Iv4}. Recall that a $1$-cycle in a complex space $Y$ is a
formal sum $Z=\sum_jn_jZ_j$, where $\{Z_j\}$ is a locally finite sequence
of irreducible analytic subsets of $Y$ of pure dimension one. The space of 
analytic $1$-cycles in $Y$ will be denoted as $\calc^{loc}_1(Y)$. It carries 
a natural topology, \ie the topology of currents. 

\smallskip From now on $Y=\Delta^2_{1+\eps}\times X$. Denote by $\calc_T$ 
the subset of $\calc^{loc}_1(Y)$ which consits of graphs $\Gamma_t$, \ie $\calc_T = 
\{\Gamma_t:t\in T\}$. We see $\calc_T$ as a topological subspace of $\calc$ and
in the sequel we shall identify $T$ with $\calc_T$. Indeed, note that 
$t\to \Gamma_t$ is injective, because such is already $t\to C_t$. Since
$T$ was supposed to be uncountable, therefore so such is also 
$\{\Gamma_t:t\in T\}=\calc_T$.

\smallskip Denote by $\bar \calc_T$ the closure of $\calc_T$ in our space of $1$-cycles 
$\calc^{loc}_1(Y)$ on $\Delta^2_{1+\eps}\times X$. Cycles $Z$ in $\bar\calc_T$ 
are characterized by following two properties:

\smallskip\sli $Z$ has an irreducible component $\Gamma$ which is a graph 
of the extension of the restriction $f|_{C\cap R_{1+\eps}}$, where $C$ is a graph
of some holomorphic function $\phi:\Delta_{1+\eps} \to \bar\Delta$.

\smallskip\slii other irreducible components fo $Z$ (if any) are a finite number
of rational curves projecting to points in $\Delta\times \Delta_{1+\eps}$.

\smallskip This directly follows from the theorem of Bishop, because areas of graphs $\Gamma_t$
are uniformly bounded, and from Lemma 7 in \cite{Iv1}, which says that a limit of a 
sequence of disks is a disk plus a finite number of rational curves. More presicely in 
(\sli we mean that $C$ is a graph of some holomorphic $\phi :\Delta_{1+\eps}\to\bar\Delta$
and $f|_{C\cap R_{1+\eps}}$ holomorphically extends to $C$ with $\Gamma_{f|_{C}}
= \Gamma$, see Lemma 1.3 from \cite{Iv4} for more details. Remark that by the choice we 
made we have that 
\begin{equation}
\eqqno(area2)
(k-1)\frac{\nu}{2} \le \area (Z)\le k\frac{\nu}{2}
\end{equation}
for all $Z\in \bar\calc_T$. Indeed \eqqref(area2) is satisfied
for $Z=\Gamma_t$ and therefore for their limits.

\begin{rema} \rm 
Let us turn attention of the reader that we write $Z$ both for an $1$-cycle as an 
analytic subset of $Y=\Delta^2_{1+\eps}\times X$ and for a corresponding {\slsf point} 
in the cycle space $\calc^{loc}_1(Y)$.
\end{rema}

\begin{figure}[h]
\centering
\includegraphics[width=3in]{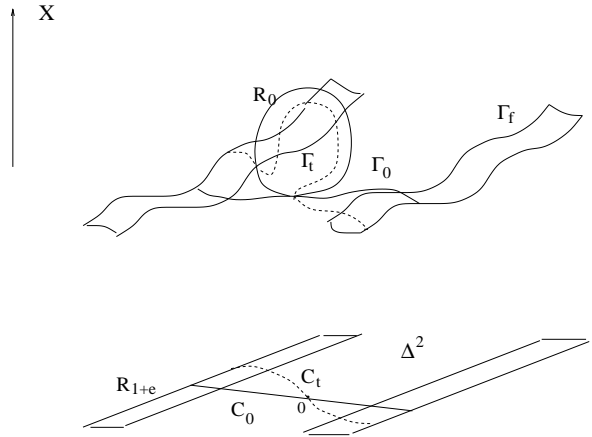}
\caption{When $C_t$ (the punctured curve downstairs) approaching $C_0$ (the bold line)
the  graph $\Gamma_t=Z_t$ (the punctured curve upstairs) of $f|_{C_t}$ stays irreducible and 
approaches $Z_0$ (the bold curve upstairs). This last is reducible, its irreducible 
component $\Gamma_0$ is a graph  over $C_0$. Its second irreducible component $R_0$ is a 
rational curve, which is contained in $\{0\}\times X$. $\Gamma_0=\Gamma_0'\cup R_0$ is 
an element of $\bar\calc_T\setminus \calc_T = \calr_T$. $\Gamma_f$ is the graph of $f$ 
over $R_{1+\eps}$.} 
\label{pic2}
\end{figure}

\smallskip Denote by  $\calr_T$ the subset  of reducible cycles in $\bar\calc_T$. This
is a closed subset of $\bar\calc_T$. Indedd, if $Z_n$ is a converging sequence from $\calr_T$ 
then every $Z_n$ has at least one irreducible component, say $R_n$, which is a rational curve.
Therefore the limit $Z\deff \lim Z_n$ contains a limit of $R\deff\lim R_n$ (up to taking a 
subsequence, if necessary). This $R$ can be only a union of rational curves. I.e., $Z$ is
reducible. The difference  $\bar\calc_T\setminus R_{\bar T}$ is uncountable since it contains 
$\calc_T$.

From here we get easily that there exists a point $Z_0\in \bar\calc_T
\setminus \calr_T$ having a fundamental system of neighborhoods $\{U_n\}$ 
in $\calc^{loc}_1(Y)$ such that $U_n\cap \bar\calc_T \subset \bar\calc_T \setminus \calr_T$ 
for all $n$ and such that all these intersections
are uncountable and relatively compact. The last is again by the theorem of Bishop.

\smallskip
Remark now that for every $Z_1, Z_2\in \bar\calc_T\cap U_1$ we have
\begin{equation}
\eqqno(area3)
\left|\area \left(Z_1\right) - \area\left(Z_2\right)
\right|\ \le  \frac{\nu}{2}.
\end{equation}
This readily follows from \eqqref(area2). First step in the proof of 
Lemma 2.4.1 from \cite{Iv3} states that a 
family of cycles satisfying \eqqref(area3) is continuous in the cycle space topology.
Remark that since $Z_1$ and $Z_2$ are irreducible we have that $Z_1= \Gamma_{f|_{C_1}}$ 
and $Z_2=\Gamma_{f|_{C_2}}$ for some discs $C_i = \{z=\phi_(\lambda)\}$, and \eqqref(area3)
means that $f|_{C_1}$ is close to $f|_{C_2}$.

\newprg[GEN.red]{The proof}

For every radius $r$ close to $1$ consider the subfamily $T_r$ of
$T$ such that $\forall t\in T_r$ function $\phi_t-\phi_{t_0}$
doesn't vanish on $\d\Delta_r$. Take a sequence $r_m\nearrow
1-\eps/2$. Suppose that $T_{r_m}$ is a most countable for every
$r_m$.  Then for all $t\in T\setminus (T_{r_1}\cup ...\cup
T_{r_m})$, which is an uncountable set, the function
$\phi_t-\phi_{t_0}$ vanishes on all $\d\Delta_{r_i}, i=1,..,m$. For
$m>N$ we get a contradiction ($N$ here bounds the winding numbers of
our test family, see discussion before Corollary \ref{levi-gen} in
Introduction). Therefore for some radius $r\sim 1-\eps/2$ (we can
suppose that $r=1$ after all) we find an uncountable subfamily
$T'\subset T$ such that for every $t\in T'$ function
$\phi_t-\phi_{t_0}$ doesn't vanish on $\d\Delta$. Take this $T'$ as
$T$ and make the reductions of the previous subsection for this $T$.

\smallskip Now remark that $Z_0$ is irreducible, \ie is an analytic disk, simply
because $Z_0$ was taken from $\bar\calc_T\setminus \calr_T$. According to (\sli 
$Z_0$ is a graph of the extension of $f|_{C_0\cap R_ {1+\eps}}$ to $C_0$ for
some curve $C_0 = \{z=\phi_0(\lambda)\}$. Take a Stein neighborhood $W$ of 
the disk $Z_0 = \Gamma_{f|_{C_0}}$, see \cite{Si2} and remark that by continuity 
of the family $\{Z: Z\in \bar\calc_T\cap U_1\}$ we have that $Z\subset W$ 
for all $Z\in U\cap (\bar\calc_T\setminus\calr_T)$ for some neighborhood 
$U\subset U_1$ of $Z_0$ in the space of cycles. Every such $Z$ is the graph 
of the extension to some $C=\{z=\phi (\lambda)\}$ of the resriction 
$f|_{C\cap R_{1+\eps}}$. Via an imbedding of $W$ to an appropriate $\cc^n$,
our $f$ is an $n$-couple holomorphic functions which holomorphically extend to 
every corresponding $C$. We are in position to apply the (holomorphic functions 
case of) Theorem \ref{levi-nonlin} and get a holomorphic
extension of $f$ to an appropriate pinched domain. This finishes the proof.

\smallskip\qed

\newprg[GEN.pos]{General position and further assumptions}

In practice one looks for extending $f$ to a bidisk $\Delta^2$. As we had seen this depends first of all
on whether a test sequence/family is in general position.  The last can be expressed in several different
ways. One of hem was given in Introduction. Another one was given in \cite{Dh} and used also in \cite{Sk}.
It sounds as follows: a family (or, a sequence) $\{C_t\}$ is said to be in general position if for any
$t_1\not=t_1\not=t_3$ one has
\begin{equation}
\eqqno(gen-pos2)
C_{t_1}\cap C_{t_2}\cap C_{t_3} = \emptyset,
\end{equation}
\ie if no three of our curves pass through one point. When $C_t$
{\slsf ad hoc} belong to a finite dimensional analytic family this
notion is equivalent to ours, simply because the set of $\alpha$
such that $f$ extends along $C_{\alpha}$ is an analytic set and
{\slsf a fortiori} forms a pinched domain, to which all  but finite
of $C_t$ should belong. In general these notions seem to be
different. Hoverer let us remark that for an uncountable family
condition \eqqref(gen-pos2) implies ours. Indeed, given
$\lambda_0\in \Delta$ if for every bidisk
$\Delta^2((\lambda_0,0),\frac{1}{n})$ the set of $t\in T$ such that
$C_t\cap \Delta^2((\lambda_0,0),\frac{1}{n}) =\emptyset$ is at most
countable then $T$ would be at most countable, unless almost all
$C_t$ pass through $(\lambda_0,0)$. Since this is forbidden by
\eqqref(gen-pos2) we see that there exists an uncountable $T'\subset
T$ such that $ C_t\cap \Delta^2((\lambda_0,0),\frac{1}{n})
=\emptyset$ for $t\in T'$. Taking a convergent sequence from $T'$ we
have that zeroes of this sequence do not accumulate to
$\lambda_0$. Applying Theorem \ref{levi-nonlin} we extend $f$ to a
pinched domain which has no pinch at $\lambda_0$. Repeating this
argument a finite number of times we extend $f$ to a neighborhood of
$\Delta\times \{0\}$.

\smallskip One can try to define the general position condition as 
such that it insures the ``non-pinching''.
Again if $\phi_k$ a priori belong to a finite dimensional family
this condition will be equivalent to the both just discussed.
Indeed, after all we know that the set of $\phi$-s such that $f$
extends along its graph is an analytic set in a finite dimensional
parameter space.
Therefore all $\phi_k$ except {\slsf finitely many} fit into a
positive dimensional families, \ie all (except finitely many) pass
(or not) through some fixed number of points. When $\phi_k$ do not
belong to a finite dimensional family (but is a test sequence) the
situation is unclear. It may happen that the Banach analytic family
$\{\phi_{\alpha}\}_{\alpha\in\cala}$ of those $\phi_{\alpha}$ along
which $f$ extend doesn't contain any of $\phi_k$. And therefore it
is not clear how to ``read off'' the ``non-pinching`` property of
the family $\{\phi_{\alpha}\}_{\alpha\in\cala}$ from the behavior of
$\phi_k$.

\smallskip For  the last point suppose now that our sequence/family 
is in general position, as in Introduction, and therefore $f$ extends
to a neighborhood of $\Delta\times \{0\}$ (or to a neighborhood of the 
graph $C_{\phi_0}$, but this is the same). The extendability of $f$ 
further to the whole of $\Delta^2$ depends now on the image space $X$.
More precisely it depends on the fact wether a Hartogs type extension 
theorem is valid for meromorphic mappings with values in this particular 
$X$. If $X$ is projective or, more generally K\"ahler, then this is true 
and was proved in \cite{Iv2}. For more general $X$ this is not always the 
case, see \cite{Iv4} for examples and further statements on this subject.

\newsect[EXMP]{Examples}

\newprg[EXMP.exmp1]{Construction of the Example 1}

\begin{nnexmp}
\label{example1}
Let the function $f$ be defined by the following series
\begin{equation}
\eqqno(exmpl1)
\sum_{n=1}^{\infty}3^{-4n^3}
\prod_{j=1}^n\left[z - \left(\frac{2}{3}\lambda \right)^j \right]
\lambda^{-n^2}z^n.
\end{equation}
Then $f$ is holomorphic in the ring domain $R\deff\cc^*\times \cc$,
holomorphically extends along every  $C_k\deff \{z=
\left(\frac{2}{3}\lambda \right)^k\}$, but there doesn't exist an
analytic family $\{\phi_{\alpha}\}_{\alpha\in \cala}$ parameterized
by a disk $\cala \ni 0$, $\phi_{0} \equiv 0$, such that
$f|_{C_{\alpha}\cap (\cc^*\times\cc)}$ meromorphically extends to
$C_{\alpha}$ for all $\alpha \in \cala$.
\end{nnexmp}
First of all the terms of this series are holomorphic and converge normally
to a holomorphic function
in the ring domain  $R = \cc^*\times \cc$. Indeed, fix any $0< \eps < 1/3$, then
for $\eps <|\lambda|<\frac{1}{\eps}$ and $|z|<\frac{1}{3\eps}$ one has

\[
\prod_{j=1}^n\left|z - \left(\frac{2}{3}\lambda \right)^j \right| \le
\prod_{j=1}^n\left(\frac{1}{\eps}\right)^j = \left(\frac{1}{\eps}
\right)^{\frac{n(n+1)}{2}},
\]
and therefore
\[
\sum_{n=1}^{\infty}3^{-4n^3}\prod_{j=1}^n\left|z - \left(\frac{2}{3}
\lambda \right)^j \right|\cdot \frac{|z|^n}{|\lambda |^{n^2}} \le
\sum_{n=1}^{\infty}3^{-4n^3-n}
\left(\frac{1}{\eps}\right)^{\frac{n(n+1)}{2}}\left(\frac{1}{\eps}
\right)^{n^2 + n} \le
\]
\[
\le
\sum_{n=1}^{\infty}3^{-4n^3-n} \left(\frac{1}{\eps}\right)^{\frac{3}{2}(n^2 + n)}.
\]
I.e., the series \eqqref(exmpl1) normally converge on compacts in
$R$ to a holomorphic function, which will be still denoted as $f(\lambda,z)$. About
$f$ let us remark that for
\begin{equation}
\eqqno(exmpl2)
z=\phi_l(\lambda) = \left(\frac{2}{3}\right)^l\lambda^l, \qquad l\ge 2,
\end{equation}
the sum in \eqqref(exmpl1) is finite and is equal to
\[
\sum_{n=1}^{l-1}3^{-4n^3}\prod_{j=1}^n\left[z -
\left(\frac{2}{3}\lambda \right)^j \right]\cdot
\frac{z^n}{\lambda^{n^2}} =
\sum_{n=1}^{l-1}3^{-4n^3}\prod_{j=1}^n\left[\left(\frac{2}{3}\lambda\right)^l
- \left(\frac{2}{3}\lambda \right)^j \right]\cdot
\left(\frac{2}{3}\right)^{nl}\lambda^{n(l-n)},
\]
with all terms being polynomials, because $l>n$ there.

\begin{prop}
\label{exmpl4} There doesn't exist a complex one-parameter analytic
family $\{\phi_{\alpha}\}_{\alpha\in \Delta}$ of holomorphic functions
in $\Delta_2$ with values in $\bar\Delta$ with $\phi_0\equiv
0$ and such that for every $\alpha\in\Delta$ the restriction
$f(\lambda, \phi_{\alpha}(\lambda))$ extends from $\Delta_2^*$ to a
meromorphic function in $\Delta_2$.
\end{prop}
\proof Suppose such family exists and let $\calp$ be a corresponding
pinched domain. All pinches of $\calp$ except at zero can be removed 
using graphs $C_k$ and Theorem \ref{levi-nonlin}. For this it is
sufficient to remark that on a small disk $\Delta_{\delta}$ around 
such pinch $\phi_k$ never vanishes and therefore our sequence is test 
on $\Delta_{\delta}$. After that by Proposition \ref{pinch0} one can 
take as our one-parameter family the family 
\begin{equation}
\eqqno(taylor1) 
\phi_{\alpha}(\lambda) = \alpha\lambda^{n_0-1}
\end{equation}
with some $n_0\ge 1$. From \eqqref(taylor1) we see that for $\lambda$ 
close to zero the image of $\phi_{\alpha}(\lambda)$ as a function of $\alpha$
will contain a disk of radius $\sim c|\lambda|^{n_0}$. Therefore for
every $\lambda\in \rr^+$ close to zero there exists $\alpha \in
\Delta_{1/2}$ such that $\phi_{\alpha}( \lambda )\in \rr^+$ and
$\phi_{\alpha}(\lambda ) \ge c\lambda^{n_0}$ for some constant
$c>0$.

\smallskip Take some $n_1>n_0$ such that $\left(\frac{2}{3}\right)^{n_1}
<\frac{c}{2}$. First of all represent our function as

\begin{equation}
\eqqno(exmpl3)
f(\lambda , z) = f_1(\lambda , z) + \prod_{j=1}^{n_1}\left[z - \left(\frac{2}{3}\lambda
\right)^j \right]f_2(\lambda , z),
\end{equation}
where 
\[
f_1(\lambda ,z) = 
\sum_{n=1}^{n_1}3^{-4n^3}\prod_{j=1}^n\left[z -
\left(\frac{2}{3} \lambda \right)^j \right]
\lambda^{-n^2}z^n 
\]
and 
\[
f_2(\lambda ,z) = \prod_{j=1}^{n_1}\left[z - \left(\frac{2}{3}\lambda
\right)^j \right] \cdot\sum_{n=n_1+1}^{\infty}3^{-4n^3}
\prod_{j=n_1+1}^n\left[z - \left(\frac{2}{3}\lambda \right)^j
\right]\lambda^{-n^2}z^n.
\]
Since $f_1$ is a rational function its restriction $f_1(\lambda ,
\phi_{\alpha}(\lambda))$ will be meromorphic in $\Delta_2$.
Therefore would  $f(\lambda , \phi_{\alpha}(\lambda))$ be
meromorphic in $\Delta_2$ we would conclude that  $f_2(\lambda ,
\phi_{\alpha}(\lambda))$ is meromorphic in $\Delta_2$ to, unless
\[
\prod_{j=1}^{n_1}\left[\phi_{\alpha}(\lambda) - \left(\frac{2}{3}\lambda \right)^j\right]
\]
is identically zero. The latter is possible only if $\phi_{\alpha}$
is one of $\phi_l$ in \eqqref(exmpl2). This is not the case and actually
by Proposition \ref{wind-fin-dim} any  complex one-parameter  family
cannot contain a converging sequence with infinitely growing winding
numbers. Therefore we have that $\phi_{\alpha}$ is not one of $\phi_l$
for all non-zero $\alpha$ small enough. Therefore $f_2(\lambda ,
\phi_{\alpha}(\lambda))$ should be meromorphic in $\Delta_2$ with
pole only at zero if  we suppose that $f(\lambda ,
\phi_{\alpha}(\lambda))$ is such. This implies that $f_2(\lambda ,
\phi_{\alpha}(\lambda))$ should be meromorphic in $\Delta_2\times
\Delta$ as function of two variables $(\lambda, \alpha)$.
But the series
\begin{equation}
\eqqno(exmpl4)
f_2(\lambda ,
\phi_{\alpha}(\lambda)) = \sum_{n=n_1+1}^{\infty}3^{-4n^3}
\prod_{j=n_1+1}^n\left[\phi_{\alpha}(\lambda) - \left(\frac{2}{3}\lambda \right)^j\right]
\lambda^{-n^2}\phi_{\alpha}(\lambda)^n
\end{equation}
representing $f_2(\lambda , \phi_{\alpha}(\lambda))$ at point
$(\lambda , \phi_{\alpha}(\lambda))\in \rr^+\times \rr^+$ can be
estimated as follows. Since
\[
\prod_{j=n_1+1}^n\left[\phi_{\alpha}(\lambda) -
\left(\frac{2}{3}\lambda \right)^j\right]\ge
\lambda^{n_0(n-n_1)}\prod_{j=n_1+1}^n\left[c -
\frac{c}{2}\lambda^{j-n_0}\right] \ge
\lambda^{n_0(n-n_1)}\left(\frac{c}{2}\right)^{n-n_1},
\]
we get that
\begin{equation}
\eqqno(exmpl5)
\sum_{n=n_1+1}^{\infty}3^{-4n^3}
\prod_{j=n_1+1}^n\left[\phi_{\alpha}(\lambda) - \left(\frac{2}{3}\lambda \right)^j\right]
\lambda^{-n^2}\phi_{\alpha}(\lambda)^n \ge
\end{equation}
\[
\sum_{n=n_1+1}^{\infty}3^{-4n^3}
\lambda^{n_0(n-n_1)-n^2}\left(\frac{c}{2}\right)^nc^n\lambda^{n_0n}
=
\sum_{n=n_1+1}^{\infty}2^{-n}3^{-4n^3}\lambda^{n_0(2n-n_1)-n^2}c^{2n}.
\]
The right hand side in \eqqref(exmpl5) grows faster than any
polynomial of $\frac{1}{\lambda}$ as $\lambda \to 0, \lambda \in \rr^+$. Therefore
$f_2(\lambda , \phi_{\alpha}(\lambda))$ has essential singularity at
$\{\lambda = 0\}$. Contradiction.

\smallskip\qed

\newprg[EXMP.exmp2]{One more example}

The following example can be found in \cite{Si1}, see p. 16.

\begin{nnexmp}\rm
\label{example2} Let $\{z_k\}_{k=0}^{\infty}$ be a sequence converging
to zero, $z_k\not= 0$. Let $P_l(z)$ be a polynomial of
degree $l+1$ such that $P_l(z_0) = ...= P_l(z_l)=0$ and $P_l(0)\not=
0$ with $\norm{P_l}_{L^{\infty}(\Delta)}= \frac{1}{l!}$. Set
\begin{equation}
\eqqno(exmple2)
f(\lambda ,z) = \sum_{l=1}^{\infty} P_l(z)\lambda^{-l}.
\end{equation}
Function $f$ is holomorphic in $\cc^*\times \cc$ and $\{0\}\times \cc$ is its
essential singularity. For every $z_k$ the restriction $f|_{C_k}\deff f(\cdot,z_k)$,
where $C_k\deff \Delta\times \{z_k\}$, is rational, having a pole of order $k$ at zero.
Moreover disks $C_k$ are test and in general position. Therefore the conclusion of Theorem
\ref{levi-nonlin} and of Corollary \ref{levi-c} fails when the orders of poles of 
restrictions are not uniformly bounded.
\end{nnexmp}

\ifx\undefined\bysame
\newcommand{\bysame}{\leavevmode\hbox to3em{\hrulefill}\,}
\fi

\def\entry#1#2#3#4\par{\bibitem[#1]{#1}
{\textsc{#2 }}{\sl{#3} }#4\par\vskip2pt}

\end{document}